# Multi-component Maintenance Optimization: A Stochastic Programming Approach

Zhicheng Zhu[1], Yisha Xiang[1], Bo Zeng[2]

[1]Department of Industrial, Manufacturing & Systems Engineering, Texas Tech University, Lubbock, TX 79409

[2]Department of Industrial Engineering, University of Pittsburgh, Pittsburgh, PA 15261

**Abstract**

Maintenance optimization has been extensively studied in the past decades. However, most of the existing maintenance models focus on single-component systems and are not applicable for complex systems consisting of multiple components, due to various interactions between the components. Multi-component maintenance optimization problem, which joins the stochastic processes regarding the failures of the components with the combinatorial problems regarding the grouping of maintenance activities, is challenging in both modeling and solution techniques, and has remained as an open issue in the literature. In this paper, we study the multi-component maintenance problem over a finite planning horizon and formulate the problem as a multi-stage stochastic integer program with decision-dependent uncertainty. There is a lack of general efficient methods to solve this type of problem. To address this challenge, we use an alternative approach to model the underlying failure process and develop a novel two-stage model without decision-dependent uncertainty. Structural properties of the two-stage problem are investigated, and a progressive-hedging-based heuristic is developed based on the structural properties. Our heuristic algorithm demonstrates a significantly improved capacity in handling practically large-size two-stage problems comparing to three conventional methods for stochastic integer programming, and solving the two-stage model by our heuristic in a rolling horizon provides a good approximation of the multi-stage problem. The heuristic is further benchmarked with a dynamic programming approach commonly adopted in the literature. Numerical results show that our heuristic can lead to significant cost savings compared with the benchmark approach.

1. **Introduction**

Effective maintenance plays an important role in maintaining high levels of productivity and safety in many capital-intensive industries, especially those that operate complex, hazardous systems, such as offshore oil and gas drilling systems, nuclear power plants, petrochemical plants, and space transport systems (Alkhamis and Yellen 1995, Cowing et al. 2004, Laggoune et al. 2009). A number of catastrophic failures, e.g., the space shuttle Challenger accident and the loss of Piper Alpha oil platform (Cowing et al. 2004), have occurred in part because of inadequate maintenance. Moreover, the downtime cost caused by either planned or unplanned maintenance shutdown in these industries is often significant. The production losses can range from $5,000 to $100,000 per hour during the shutdown in chemical plants and millions of dollars per day in offshore drilling/refineries (Tan and Kramer 1997, Amaran et al. 2016). As the demand for high reliability increases, it is more imperative to develop efficient maintenance schedules for complex systems.

Maintenance optimization has been extensively studied in the literature. However, most of the existing maintenance models focus on single-component systems, and are not applicable for complex systems consisting of multiple components, due to various interactions between the components. In general, there are three different types of interactions, economic, structural, and stochastic dependence (Thomas 1986). Economic dependence is the most common one among these three types of interactions. Systems with the economic dependence typically incur a common system-level cost, often referred to as setup cost, due to mobilizing repair crew, safety provisions, disassembling machines, special transportation, and the downtime loss. These costs are shared by all maintenance activities performed simultaneously. Considerable cost savings can be obtained by jointly maintaining several components instead of separately, especially when the setup cost is high.

Multi-component maintenance optimization problem, which joins the stochastic processes regarding the failures of the components with the combinatorial problems regarding the grouping of maintenance activities (Dekker et al. 1997, Scarf 1997, Dekker and Scarf 1998, Van Horenbeek and Pintelon 2013), is challenging in both modeling and solution techniques, and has remained as an open issue in the literature. The problem can quickly end in complex models and explicit analytical expressions for optimal maintenance costs and the corresponding decisions are sometimes impossible to obtain. One often has to make special system assumptions (Castanier et al. 2005, Tian and Liao 2011, Huynh et al. 2015), impose restrictions on maintenance grouping activities (Dekker et al. 1997, Wildeman and Dekker 1997, Ding and Tian 2012, Van Horenbeek and Pintelon 2013), and/or resort to simulation tools (Tan and Kramer 1997, Bérenguer et al. 2000, Barata et al. 2002, Laggoune et al. 2009, Nguyen et al. 2015) so that the decision problem can be formulated with less mathematical difficulty.

In this paper, we study the maintenance optimization problem for multi-component systems with economic dependence over a finite planning horizon. We aim to minimize total maintenance cost by optimally determining maintenance decisions for all components at each decision period. We first formulate this problem as a multi-stage stochastic integer program. The system state transition probabilities at each stage are determined by not only the underlying failure processes but also maintenance decisions, and thus are decision-dependent. Such decision-dependent uncertainty is also referred to as endogenous uncertainty, and there is a lack of efficient methods to handle this type of problem (Peeta et al. 2010, Zhan et al. 2016, Apap and Grossmann 2017). We approximate the multi-stage model with a novel two-stage stochastic linear integer model in a rolling horizon. In both models, we do not restrict the types of maintenance activities that can be grouped or when the grouping can occur. In other words, joint execution of any combination of maintenance activities can occur anytime. A progressive-hedging-based heuristic is designed to solve practically large-size two-stage problems. Computational studies are performed to assess the performance of the heuristic. We further compare the results of the two-stage model with those of a direct-grouping model using a dynamic-programming approach (Van Horenbeek and Pintelon 2013) over the rolling horizon. The main contributions of this paper are as follows.

(1) From a modeling perspective, the proposed multi-stage model and the two-stage model are sufficiently general, permitting grouping of any maintenance activities at any time and allowing general failure distributions. This work extends the multi-component maintenance literature by using a stochastic programming approach. Stochastic programming is a powerful modeling technique and facilitates the derivation of analytical expressions of the total cost function and maintenance decisions, which are difficult to obtain using commonly adopted approaches for the multi-component problem, such as dynamic programming.

(2) Formulating the multi-component maintenance problem as a stochastic program further enables the use of stochastic optimization tools. We design an efficient heuristic algorithm under the progressive hedging framework based on the problem structural properties of the two-stage model. Our heuristic algorithm demonstrates a significantly improved capacity in handling practically large-size two-stage problems comparing to three conventional methods for stochastic integer programming, and solving the two-stage model by our heuristic in a rolling horizon provides a good approximation of the multi-stage problem. Using a rolling horizon scheme, we further assess the performance of our heuristic by comparing it with a dynamic programming approach widely adopted in the literature. Numerical results show that the two-stage maintenance model and the designed heuristic can lead to substantial cost savings.

The remainder of this paper is organized as follows. Section 2 reviews multi-component maintenance and important solution methods in stochastic programming. In Section 3, we develop a multi-stage model stochastic maintenance model and approximate it with a two-stage model. Section 4 describes the progressive-hedging-based heuristic and three conventional algorithms in detail. Computational studies are presented in Section 5. We conclude this study and discuss future research directions in Section 6.

## 2. Literature review

We review two streams of literature that are relevant to our work: literature on multi-component maintenance and literature on solution methods in stochastic programming.

### 2.1 Multi-component maintenance

Common approaches to coordinating maintenance activities of multi-components include direct-grouping, indirect-grouping and opportunistic maintenance. Direct-grouping partitions the components into a number of fixed groups and then always maintain the components in a group jointly (van Dijkhuizen and van Harten 1996). The problem formulated with this approach is a NP-complete set-partitioning problem. The optimal grouping decision can be found for only a small number of components due to the computational complexity. There are some efforts that reduce the set-partitioning problem for multi-component maintenance to a dynamic-programming problem with a quadratic time complexity under some special assumptions (Dekker et al. 1996, Wildeman et al. 1997). Van Horenbeek and Pintelon (2013) and Vu et al. (2014) extend Dekker et al. (1996) Wildeman et al. (1997) by considering dynamic information, e.g., usage of components and environmental conditions. A major deficiency of this approach is that grouping activity iteratively takes place within a time window that is often determined by the maximum individual maintenance interval among all components. Within this window, each component is preventively maintained only one time. This assumption is not relevant since a system may be composed of different components with different lifetime cycles, and maintenance intervals of components can be significantly different (Laggoune et al. 2009).

Unlike the direct-grouping that yields a fixed group structure, indirect grouping usually groups preventive maintenance (PM) activities by making the PM interval a multiple of a basis interval, so the maintenance of different components can coincide (Goyal and Kusy 1985, Goyal and Gunasekaran 1992, Vos de Wael 1995). An alternative indirect grouping strategy performs major PM on all components jointly at the end of a common interval and allows minor or major PM within this interval. Indirect grouping model of this type is sometimes formulated as a mixed integer programing (MIP) problem (Epstein and Wilamowsky 1985, Hariga 1994). Because of the simplified policy structure, the MIP model can be separated by components, which greatly reduces the computational complexity.

Both direct- and indirect-grouping focus on grouping PM activities, and ignore maintenance opportunities generated by corrective maintenance (CM) at failures. To take advantage of the time window of CM and use it as opportunities for PM of other functioning components, many opportunistic maintenance (OM) models have been proposed (Pham and Wang 2000, Rao and Bhadury 2000, Cui and Li 2006, Besnard et al. 2009, Ding and Tian 2012, Koochaki et al. 2012, Patriksson et al. 2015). Ding and Tian (2012) and Koochaki et al. (2012) use a simulation-based optimization method to find optimal OM policies. Shafiee and Finkelstein (2015) consider a simplified OM policy that preventively replaces all non-failed components when there is a failure. More recently, Patriksson et al. (2015) use a stochastic programming approach in OM. However, the integer L-shaped method proposed in Patriksson et al. (2015) cannot solve large-scale problems. Castanier et al. (2005) consider a condition-based OM policy and formulate it as a semi-regenerative process. This approach also suffers the computational intractability, because the problem size grows exponentially as the number of components increases. As a result, their analysis is limited to a two-component system. For more details regarding the multi-component maintenance problem, the readers are referred to review papers by Thomas (1986), Dekker et al. (1997), and Nicolai and Dekker (2008).

**2.2 Solution methods in stochastic programming**

In this paper, we formulate the multi-component maintenance optimization problem as a stochastic integer program. Various decomposition methods have been developed to solve stochastic integer programs. Benders decomposition (Birge and Louveaux 2011, Bodur et al. 2016, Rahmaniani et al. 2017) and progressive hedging algorithm (PHA) (Rockafellar and Wets 1991, Watson and Woodruff 2011, Gade et al. 2016) are two important decomposition methods for solving stochastic integer programming problems. Benders decomposition vertically decomposes a problem into a master problem that only concerns first-stage decisions and subproblems that include second-stage decisions of all scenarios. Benders cut and integer L-shaped cut are two major types of cuts that are added within Benders decomposition framework. However, Benders cut may become useless since strong duality does not hold in an integer program, and integer L-shaped cut is typically inefficient because every feasible solution may need an integer L-shaped cut in the worst-case scenario. The PHA decomposes the extensive form according to scenario, and iteratively solving penalized versions of the sub-problems to gradually enforce non-anticipativity. The performance of PHA, to a great extent, depends on how efficient each subproblem is solved in a stochastic integer program. In our problem, each scenario subproblem with deterministic individuals' lifetimes is essentially a NP-complete set-partitioning problem. Efficient heuristic algorithm for each subproblem is needed.

Our review of the literature shows that few research has considered grouping at both preventive and corrective maintenance occasions under practical assumptions, which significantly affects the optimality of the solutions because of the simplified models and reduced solution space. There is also a lack of efficient algorithms that can provide satisfactory results for practically large-scale multi-component maintenance problems.

## 3. Model development

**Notation**

| | |
|---|---|
| **Parameters** | |
| $n$ | Number of components |
| $N$ | Component set, $N = \{1, 2, …, n\}$ |
| $T$ | Length of the planning horizon |
| $T_s$ | Planning horizon, $T_s = \{1, 2, …, T\}$ |
| $q$ | Number of individuals |
| $R$ | Individual set, $R = \{1, 2, …, q\}$ |
| $\tilde{\Omega}$ | Set that consists of all possible scenarios in the two-stage model |
| $\Omega$ | Scenario set that sampled from $\tilde{\Omega}$ |
| $I_{ir}$ | Individual $r$ of component $i$ |
| $T_{ir}$ | Lifetime of $I_{ij}$ |
| $T_{ir}^{\omega}$ | Lifetime of $I_{ij}$ in scenario $\omega$ |
| $T'$ | Extended planning horizon, $T' = \max_{i,r,\omega} T_{ir}^{\omega}$ |
| $c_{i,\text{pr}}$ | Preventive replacement (PR) cost for type $i$ component's individual |
| $c_{i,\text{cr}}$ | Corrective replacement (CR) cost for type $i$ component's individual |
| $C_{i,\text{pr}}$ | Total PR cost incurred by individuals of component $i$ in the planning horizon $T_s$ |
| $C_{i,\text{cr}}$ | Total CR cost incurred by individuals of component $i$ in the planning |
|  | horizon $T_s$ |
| $C_s$ | Total setup cost in the planning horizon $T_s$ |
| $d$ | Setup cost |
| $\delta$ | Length of a decision period |
| $\xi_{it}$ | Equal to 1 when the working individual of component $i$ is failed at the stage $t$, otherwise 0 |
| $\xi_t$ | Vector of all working individuals' states at stage $t$, $\xi_t = \{\xi_{1,t}, \ldots, \xi_{2,t}, \ldots, \xi_{nt}\}$ |
| $a_{it}$ | Age of the working individual of component $i$ prior to maintenance decision at stage $t$ |
| $a_t$ | Vector of all working individuals' ages at stage $t$, $a_t = \{a_{1,t}, \ldots, a_{2,t}, \ldots, a_{nt}\}$ |
| $Q(x)$ | Expected second-stage cost |
| $p(\omega)$ | Probability of scenario $\omega$ |

**Decision Variables**

|  |  |
|---|---|
| $x_i$ | Equal to 1 when an individual of component $i$ is replaced at the first stage, 0 otherwise |
| $z$ | Equal to 1 when there is at least one individual maintained at the first stage, 0 otherwise |
| $x_{it}$ | Equal to 1 when an individual of component $i$ is replaced at stage $t$, 0 otherwise |
| $z_t$ | Equal to 1 when there is at least one individual maintained at stage $t$, 0 otherwise |
| $\tilde{x}_{it}^{r\omega}$ | Equal to 1 when $I_{ir}$ is replaced at or before time $t$ in scenario $\omega$, 0 otherwise |
| $z_t^{\omega}$ | Equal to 1 when there is at least one individual maintained at time $t$ in scenario $\omega$, 0 otherwise |

## 3.1 Problem statement

In this paper, we consider maintenance optimization for a multi-component system with economic dependence. The objective is to minimize total expected maintenance cost in a finite planning horizon. We consider two types of maintenance, preventive replacement (PR) and corrective replacement (CR). At each decision period, maintenance decisions need to be made for all components. Any maintenance activities of any components can be performed together to save the setup cost and improve the system performance. Note that any failed component is correctively replaced. Both CR and PR use a new component and are therefore perfect. Components' failure models can be any failure distribution.

### 3.2 Multi-stage stochastic maintenance model

Consider a system that consists of $N = \{1, \ldots, n\}$ components. Each component in the system is considered as a different type of component regardless of its physical type. To distinguish the component type and the component itself, we use component only when referring to its type and refer to physical components as individuals. For example, individual $I_{ir}$ is the individual used for the $r^{th}$ maintenance replacement of type $i$ component.

We consider a discretized finite planning horizon $T_s = \{1, \ldots, T\}$, where the length of a decision period is $\delta$. Denote the PR cost and CR cost for each component by $c_{i,\text{pr}}$ and $c_{i,\text{cr}}$ respectively, and assume $c_{i,\text{pr}} < c_{i,\text{cr}}$ for all $i \in N$. The system setup cost is $d$ at any maintenance occasion. If $n$ individuals are replaced at the same time, the total savings from executing these $n$ maintenance activities jointly is $d(n-1)$.

This problem is naturally a multi-stage stochastic programming problem. Denote the state for the working individual of component $i \in N$ at stage $t \in T_s$ by $\xi_{it}$, i.e., $\xi_{it} = 1$ if the working individual fails and 0 otherwise. Let $a_{it}$ be the age of the working individual of component $i \in N$ prior to the maintenance decision at stage $t$, and $a_t = (a_{1,t}, a_{2,t}, \ldots, a_{nt})$ be the age vector of all working individuals. At each decision stage $t$, after observing the states $\xi_t = (\xi_{1,t}, \xi_{2,t}, \ldots, \xi_{nt})$ of all working individuals, we first correctively replace all failed individuals and then select a group of individuals for PR if desired.

Denote the decision variables by

$$x_{it} = \begin{cases} 1, & \text{if the individual of component } i \text{ is replaced at time } t, \quad i \in N, t \in T_s \\ 0, & \text{otherwise,} \quad\quad\quad\quad\quad\quad\quad\quad\quad\quad\quad\quad\quad\quad\quad\quad i \in N, t \in T_s \end{cases}$$

and

$$z_t = \begin{cases} 1, & \text{if any maintenance occurs at time } t, \ t \in T_s \\ 0, & \text{otherwise,} \hspace{3em} t \in T_s \end{cases}.$$

The multi-stage stochastic model **(P1)** is defined as follows:

$$f_1(a_1,\xi_1) = \min \sum_{i \in N} c_{i,\text{pr}} x_{i,1} + \sum_{i \in N} (c_{i,\text{cr}} - c_{i,\text{pr}})\xi_{i,1} + dz_1 + V_1(x_1,a_1) \tag{1a}$$

subject to

$$V_t(x_t,a_t) = \begin{cases} E_{\xi_{t+1}}[f_{t+1}(a_{t+1},\xi_{t+1})], & t \in T_s \setminus \{T\} \\ 0, & t = T \end{cases} \tag{1b}$$

$$a_{t+1} = a_t(1-x_t) + \delta, \quad t \in T_s \setminus \{T\} \tag{1c}$$

$$x_{it} \geq \xi_{it}, \quad i \in N, t \in T_s \tag{1d}$$

$$z_t \geq x_{it}, \quad i \in N, t \in T_s \tag{1e}$$

$$x_{it} \in \{0,1\}, \ i \in N, t \in T_s \tag{1f}$$

$$z_t \in \{0,1\}, \ t \in T_s \tag{1g}$$

Objective function (1a) includes maintenance cost at the first stage and the expected minimum cost at the second stage. The expected minimum second-stage cost is given by Constraints (1b). Constraints (1c) provide the age of each work individual at each decision stage. Constraints (1d) ensure that the indicator of replacement $x_{it}$ is 1 when an individual failed. Constraints (1e) force that setup cost is incurred whenever a replacement is performed. Constraints (1f) and (1g) are integrality constraints.

The maintenance decision at each decision stage influences the system state transition probability. Scenario tree in Figure 1 illustrates the interactions between maintenance decision and the underlying stochastic failure process of a two-component system.

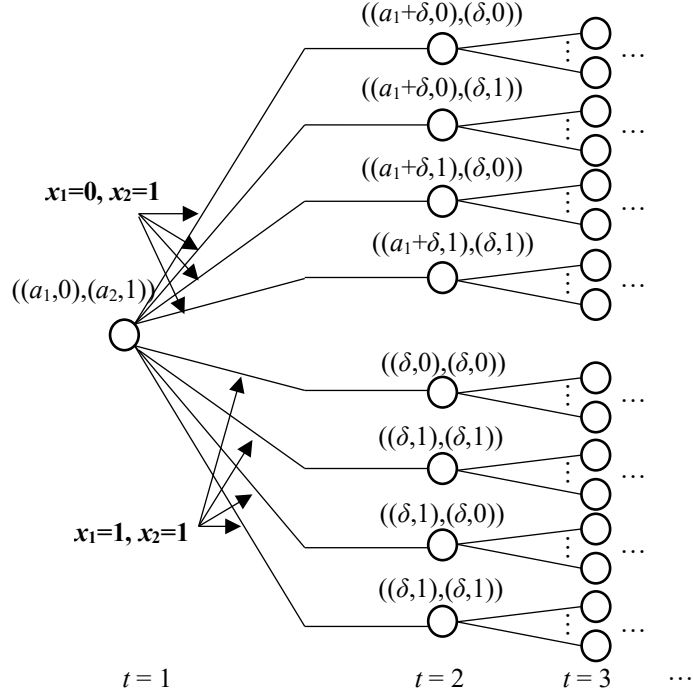

Figure 1. Scenario tree with different decisions

Model P1 is therefore a stochastic integer program with decision-dependent uncertainty, where decision-dependent uncertainty is also referred to as endogenous uncertainty. This type of problem is difficult to solve. First, general efficient method for stochastic integer programs is lacking. Second, the problem size grows exponentially as the number of components and/or decision stages increases. Third, there is a lack of general efficient algorithms to solve a stochastic program with endogenous uncertainty. Next, we use a novel two-stage model to approximate P1 and design an efficient heuristic to find high-quality solutions.

### 3.3 Two-stage stochastic maintenance model

A common approach to approximating a multi-stage stochastic program is to utilize a two-stage model in a rolling horizon. The two-stage approximation model is usually obtained by combining all future periods together as a second stage problem with all future's non-anticipativity constraints removed (Mulvey and Vladimirou 1992, Tolio and Urgo 2007, Beraldi et al. 2011).

In this two-stage model, instead of using failure probability as the underlying stochastic element to generate scenarios, we use the random lifetimes as the equivalent stochastic element to capture the uncertainty. Specifically, let $T_{ir}^{\omega}$ be the lifetime of individual $I_{ir}$ that drawn from its failure distribution, a scenario $\omega \in \tilde{\Omega}$ is a lifetime combination of all individuals of all components, i.e., $(T_{1,1}^{\omega},...,T_{1,q}^{\omega},...,T_{n,1}^{\omega},...,T_{nq}^{\omega})$, where $\tilde{\Omega}$ is the collection of all possible scenarios. Note that the endogenous uncertainty is removed by using this alternative approach to describe the uncertainty.

The two-stage model can be represented by

$$\min_{x \in X} \sum_{\omega \in \tilde{\Omega}} p(\omega) F(x, \omega),$$

where $x$ is the vector of decision variables and $X$ is the feasible region. Function $F(x, \omega)$ is the total maintenance cost given the lifetime combination in scenario $\omega \in \tilde{\Omega}$, and $p(\omega)$ is the probability of scenario $\omega$. Since the realization of lifetime $T_{ir}^{\omega}$ is infinite for the majority of lifetime distributions, e.g., Weibull distribution, the total number of scenarios $|\tilde{\Omega}|$ is infinite. We therefore use the sample average approximation (SAA) method (Kleywegt et al. 2002) to approximate the two-stage model. Specifically, we have

$$\sum_{\omega \in \tilde{\Omega}} p(\omega) F(x, \omega) \approx \sum_{\omega \in \Omega} p(\omega) F(x, \omega),$$

where $\Omega \subset \tilde{\Omega}$ and $p(\omega) = \frac{1}{|\Omega|}$.

Before presenting function $F(x, \omega)$, we first introduce decision variables in the two-stage model. Let

$$x_i = \begin{cases} 1, & \text{if an individual of component } i \text{ is replaced at the first stage,} \quad i \in N \\ 0, & \text{otherwise,} \quad i \in N \end{cases}$$

$$\tilde{x}_{it}^{r\omega} = \begin{cases} 1, & \text{if } I_{ir} \text{ is replaced at or before time } t \text{ in scenario } \omega, \quad i \in N, t \in T_s, r \in R, \omega \in \Omega \\ 0, & \text{otherwise.} \quad i \in N, t \in T_s, r \in R, \omega \in \Omega \end{cases}$$

and

$$z_t^\omega = \begin{cases} 1, & \text{if maintenance occurs at time } t \text{ in scenario } \omega, t \in T_s, \omega \in \Omega \\ 0, & \text{otherwise.} \quad t \in T_s, \omega \in \Omega \end{cases}.$$

To facilitate the model development, we introduce two auxiliary binary variables $Y_i^{r\omega}$ and $w_{it}^{r\omega}$ based on $\tilde{x}_{it}^{r\omega}$. The variable $Y_i^{r\omega}$ is an indicator of the maintenance type. More details regarding these auxiliary variables will be discussed in Section 3.3.1. The deterministic extensive form (DEF) of the two-stage model, which explicitly describes the second-stage decision variables for all scenarios (Birge and Louveaux 2011), is formulated as follows.

**Model DEF**:

minimize

$$\sum_{\omega \in \Omega} p(\omega) \left( \sum_{i \in N} \left( \underbrace{\sum_{r=1}^q \left( c_{i,\text{pr}} Y_i^{r\omega} \right)}_{C_{i,\text{pr}}} + \underbrace{\sum_{r=1}^q \left( c_{i,\text{cr}} \left(1 - Y_i^{r\omega}\right)\right) - c_{i,\text{cr}} \left(1 - \tilde{x}_{iT}^{r\omega}\right)}_{C_{i,\text{cr}}} \right) + \underbrace{\sum_{t \in T_s} d z_t^\omega}_{C_s} \right) \quad (2a)$$

subject to

$$\tilde{x}_{it}^{r\omega} \leq \tilde{x}_{i,t+1}^{r\omega}, \qquad i \in N, t \in T_s \setminus \{T\}, r \in R, \omega \in \Omega \quad (2b)$$

$$\tilde{x}_{i,t+1}^{r+1,\omega} \leq \tilde{x}_{it}^{r\omega}, \qquad i \in N, t \in T_s \setminus \{T\}, r \in R \setminus \{q\}, \omega \in \Omega \quad (2c)$$

$$\sum_{r \in R} \left( \tilde{x}_{it}^{r\omega} - \tilde{x}_{i,t-1}^{r\omega} \right) \leq z_t^\omega, \qquad i \in N, t \in T_s \setminus \{1\}, \omega \in \Omega \quad (2d)$$

$$\tilde{x}_{i,1}^{1,\omega} \leq z_1^\omega, \qquad i \in N, \omega \in \Omega \quad (2e)$$

$$\tilde{x}_{it}^{r\omega} \leq \tilde{x}_{i,t+T_{i,r+1}^\omega}^{r+1,\omega}, \qquad i \in N, t \in \{1, ..., T - T_{i,r+1}^\omega\}, r \in R \setminus \{q\}, \omega \in \Omega \quad (2f)$$

$$\tilde{x}_{i,T_{i1}^\omega}^{1,\omega} = 1, \qquad i \in \{j \in N \mid T_{j1}^\omega \leq T\}, \omega \in \Omega \quad (2g)$$

$$\tilde{x}_{i,1}^{r\omega} = 0, \qquad i \in N, r \in R \setminus \{1\}, \omega \in \Omega \quad (2h)$$

$$x_i = \tilde{x}_{i,1}^{1,\omega}, \qquad i \in N,\ \omega \in \Omega \tag{2i}$$

$$x_i \geq \xi_i, \qquad i \in N \tag{2j}$$

$$Y_i^{1,\omega} = 1 - w_{i,T_{i1}^{\omega}}^{1,\omega}, \qquad i \in N,\ \omega \in \Omega \tag{2k}$$

$$Y_i^{r\omega} = \left( \sum_{t=T_{ir}^{\omega}}^{T+T_{ir}^{\omega}} |y_{it}^{r\omega}| + \sum_{t=1}^{T_{ir}^{\omega}-1} w_{it}^{r\omega} \right) / 2, \qquad i \in N,\ r \in R \setminus \{1\},\ \omega \in \Omega \tag{2l}$$

$$y_{it}^{r\omega} = w_{it}^{r\omega} - w_{i,t-T_{ir}^{\omega}}^{r-1,\omega}, \qquad i \in N,\ r \in R \setminus \{1\},\ t \in \{T_{ir}^{\omega}, ..., T'\},\ \omega \in \Omega \tag{2m}$$

$$w_{it}^{r\omega} = \tilde{x}_{it}^{r\omega} - \tilde{x}_{i,t-1}^{r\omega}, \qquad i \in N,\ r \in R,\ t \in T_s \setminus \{1\},\ \omega \in \Omega \tag{2n}$$

$$w_{i0}^{r\omega} = \tilde{x}_{i0}^{r\omega}, \qquad i \in N,\ r \in R,\ \omega \in \Omega \tag{2o}$$

$$w_{it}^{r\omega} = 0, \qquad i \in N,\ r \in R,\ t \in \{T+1, ..., T'\},\ \omega \in \Omega \tag{2p}$$

$$\tilde{x}_{it}^{r\omega} \in \{0,1\}, \qquad i \in N,\ r \in R,\ t \in T_s,\ \omega \in \Omega \tag{2q}$$

$$x_i \in \{0,1\}, \qquad i \in N \tag{2r}$$

$$z_t^{\omega} \in \{0,1\}, \qquad t \in T_s,\ \omega \in \Omega \tag{2s}$$

$$w_{it}^{r\omega} \in \{0,1\}, \qquad i \in N,\ r \in R,\ t \in \{1, ..., T'\},\ \omega \in \Omega \tag{2t}$$

$$Y_i^{r\omega} \in \{0,1\}, \qquad i \in N,\ r \in R,\ \omega \in \Omega \tag{2u}$$

Function (2a) is the objective function. Decision variables $x_i$ and $z$ concern maintenance decisions at the first stage, and $\tilde{x}_{it}^{r\omega}$ and $z_t^{\omega}$ are the second-stage decisions for scenario $\omega$. We provide detailed derivation of the objective function in Section 3.3.1 and explain the constraints in Section 3.3.2.

### 3.3.1 Derivation of the Objective Function

In objective function (2a), the total cost includes: (1) sum of the PR and CR costs incurred by individuals of component $i$ in the planning horizon, denoted by $C_{i,\text{pr}}$ and $C_{i,\text{cr}}$ respectively, and (2) total system setup cost $C_s$. We break the derivation of the total cost function into the calculations of these cost elements.

- **Derivation of $C_{i,\text{pr}}$**

For component $i$, the total cost of individuals preventively replaced over the planning horizon is given by $C_{i,\text{pr}} = \sum_{r=1}^{q} c_{i,\text{pr}} Y_i^{r\omega}$, where $Y_i^{r\omega}$ is defined in constraints (2k) and (2l).

Next, we explain why $Y_i^r$ can be used to identify the replacement type and this determination is a key element of the model development. We drop the superscript $\omega$ in the following discussions for notational convenience. It is obvious that the decision variables $\tilde{x}_{it}^{r\omega}$ and $w_{it}^{r\omega}$ only concern when a placement is performed and have no indication on the type of replacement. For an individual $I_{ir}$, one way to determine its replacement type is to examine the time interval between the replacements of individuals $I_{i,r-1}$ and $I_{ir}$. Suppose that individuals $I_{i,r-1}$ and $I_{ir}$ are replaced at times $t_1$ and $t_2$ (i.e., $w_{it_1}^{r-1} = 1$ and $w_{it_2}^{r} = 1$), respectively. If the difference between $t_2$ and $t_1$ equals to the lifetime of $I_{ir}$, namely $T_{ir}$, then $I_{ir}$ is replaced at the end of its lifetime and the replacement type is CR. The replacement is PR otherwise. Therefore, if CR is performed on this individual, we have $w_{it}^r - w_{i,t-T_{ir}}^{r-1} = 0 \ \forall t \in T_s$, which leads to $\sum_{t=1}^{T}|y_{it}^r| = \sum_{t=1}^{T}|w_{it}^{r\omega} - w_{i,t-T_{ir}^\omega}^{r-1,\omega}| = 0$ (Figure 2(a)). If PR is performed on this individual, then $w_{i,t_2}^r - w_{i,t_2-T_{ir}}^{r-1} = 1$, $w_{i,t_1+T_{ir}}^r - w_{i,t_1}^{r-1} = -1$, and $w_{i,t+T_{ir}}^r - w_{it}^{r-1} = 0$ for all $t \in \{t \mid t \in T_s, t \neq t_1, t \neq t_2\}$ (Figure 2(b)), and consequently, $\sum_{t=1}^{T}|y_{it}^r| = 2$. This makes the value of $\sum_{t=1}^{T}|y_{it}^r|/2$ a good indicator for determining the replacement type, and $\sum_{t=1}^{T}|y_{it}^r|$ is calculated as follows:

$$\sum_{t=1}^{T}|y_{it}^r| = \sum_{t=1}^{T}|w_{it}^r - w_{i,t-T_{ir}}^{r-1}| \tag{3}$$

However, there are two boundary issues in Equation (3).

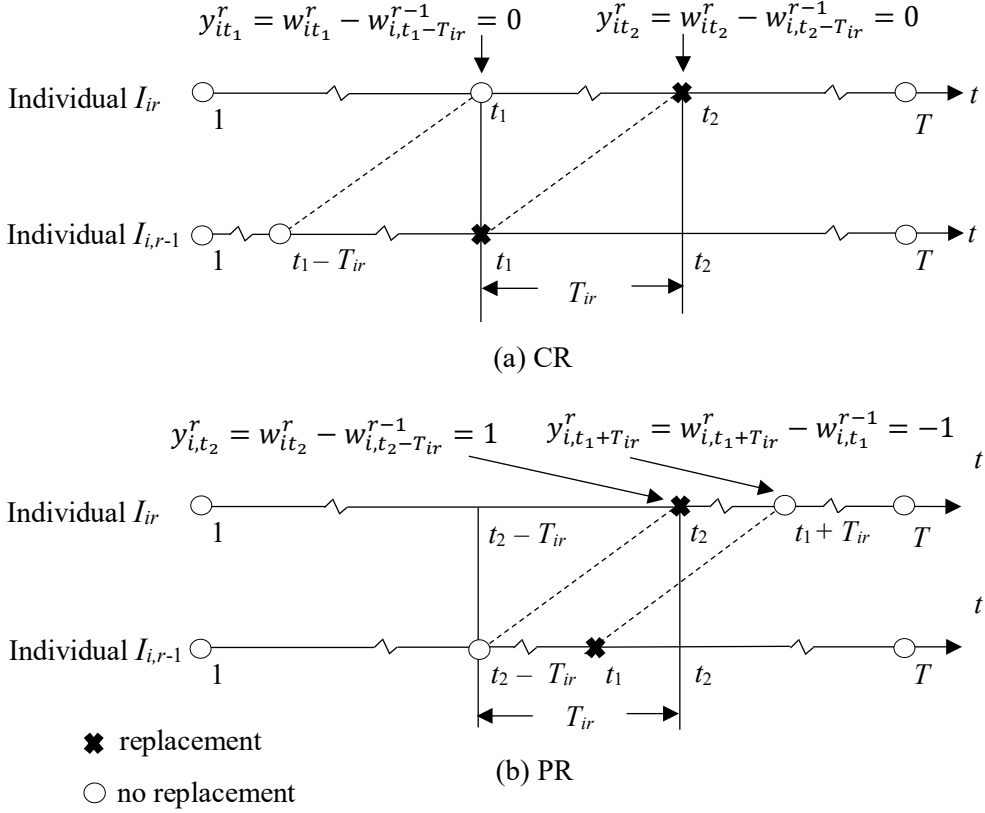

Figure 2: Illustration of distinguishing PR and CR

(1) Decision times of $w_{it}^r$ need to be extended beyond $T_s$. This is because Equation (3) does not count the maintenance decisions made for individual $I_{i,r-1}$ at times $\{T - T_{ir} + 1, \ldots, T\}$. To include these decisions, the planning horizon for $w_{it}^r$ is extended to $T' = T + \max_{i,r} T_{ir}$, and let $w_{it}^r = 0$ for $t > T$. See the region labeled "Not Defined" in Figure 3 for an illustration.

(2) In Equation (3), the decision times considered for individual $I_{ir}$ implicitly start from $T_{ir}$, and all decisions made before $T_{ir}$ are excluded (illustrated in the region labeled "Excluded"). To recover decisions made for individual $I_{ir}$ at times $\{1, \ldots, T_{ir} - 1\}$, we add $\sum_{t=1}^{T_{ir}-1} w_{it}^r$ to Equation (3). Equation (3) is now rewritten as follows,

$$Y_i^{r\omega} = \sum_{t=T_{ir}^\omega}^{T+T_{ir}^\omega} |y_{it}^{r\omega}| + \sum_{t=1}^{T_{ir}^\omega - 1} w_{it}^{r\omega}, \quad i \in N, r \in R \setminus \{1\}, \omega \in \Omega.$$

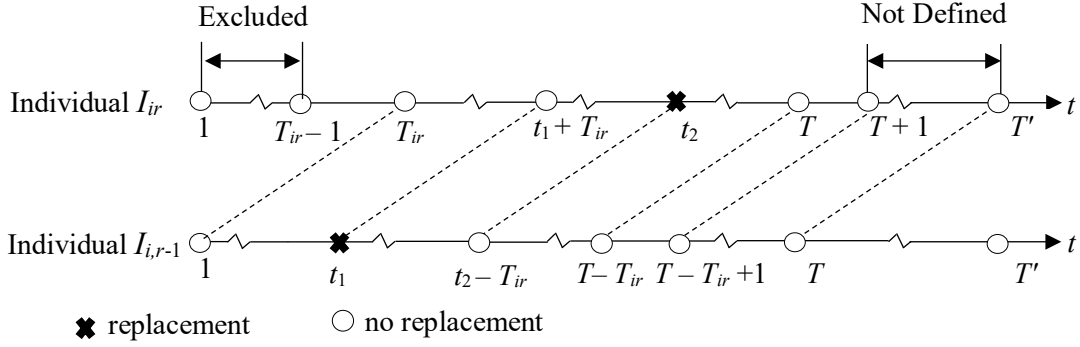

Figure 3: Illustration of the boundary issue of Equation (3)

The absolute function, $\left|y_{it}^{r\omega}\right|$, can be linearized by a pair of deviation variables $u_{it}^{r\omega}$ and $v_{it}^{r\omega}$ (Rardin and Rardin 2016). We replace $\left|y_{it}^{r\omega}\right|$ with Equation (4) in the constraint (2l), and add constraints (2v) to (2x) in the DEF model. Note that constraint (2w) is not needed for linearization but it makes the problem formulation stronger when binary integer restrictions on $u_{it}^{r\omega}$ and $v_{it}^{r\omega}$ are relaxed.

$$\left|y_{it}^{r\omega}\right| = u_{it}^{r\omega} + v_{it}^{r\omega}, \qquad i \in N,\ r \in R,\ t \in \{1,...,T'\},\ \omega \in \Omega \qquad (4)$$

$$y_{it}^{r\omega} = u_{it}^{r\omega} - v_{it}^{r\omega}, \qquad i \in N,\ r \in R,\ t \in \{1,...,T'\},\ \omega \in \Omega \qquad (2v)$$

$$u_{it}^{r\omega} + v_{it}^{r\omega} \leq 1, \qquad i \in N,\ r \in R,\ t \in \{1,...,T'\},\ \omega \in \Omega \qquad (2w)$$

$$u_{it}^{r\omega}, v_{it}^{r\omega} \in \{0,1\}, \qquad i \in N,\ r \in R,\ t \in \{1,...,T'\},\ \omega \in \Omega \qquad (2x)$$

- **Derivation of $C_{i,\mathrm{cr}}$ and $C_\mathrm{s}$**

For component $i$, the total cost of individuals correctively replaced over the planning horizon is given by $C_{i,\mathrm{cr}} = \sum_{r=1}^{q}\left(c_{i,\mathrm{cr}}\left(1-Y_i^{r\omega}\right) - c_{i,\mathrm{cr}}\left(1-\tilde{x}_{iT}^{r\omega}\right)\right)$. As explained in the derivation of $C_{i,\mathrm{pr}}$, the expression $Y_i^r = 0$ implies a CR for individual $I_{ir}$. However, for any component type, the number of individuals used for replacement is unknown due to the unknown maintenance decisions, and the maximum number of individuals needed for any component is considered in the optimization model. It is likely that some individuals are not used in the planning horizon. If neither individual $I_{i,r-1}$ nor $I_{ir}$ is used for replacement during the planning horizon, the value of $Y_i^r$ is also zero. We need to distinguish these two scenarios that both have $Y_i^r = 0$. This can be done by examining the value of $\tilde{x}_{iT}^{r\omega}$. If an individual is not used, we have $\tilde{x}_{iT}^{r\omega} = 0$ and $\tilde{x}_{iT}^{r\omega} = 1$ otherwise. The false corrective cost caused by an individual that is not used is $c_{i,\text{cr}}(1 - \tilde{x}_{iT}^{r\omega})$, and needs to be subtracted from the total cost, $C_{i,\text{cr}}$. Lastly, the total setup cost over the planning horizon is $C_s = \sum_{t \in T_s} dz_t^\omega$.

### 3.3.2 Constraints

Constraints (2a) are the definition of $\tilde{x}_{it}^{r\omega}$, which ensure that individual $I_{ir}$ is replaced at or before $t + 1$ when it is replaced at or before $t$. Constraints (2c) imply that individual $I_{i,r+1}$ can only be replaced after $I_{ir}$ is replaced. Constraints (2d) and (2e) ensure that maintenance cost $d$ incurs when any component is replaced at time $t$. Constraints (2f) and (2g) ensure that individual $I_{ir}$ is replaced at the latest when it has been inside the system for $T_{ir}^\omega$ time units. In other words, $I_{ir}$ has to be replaced before or at the end of its lifetime. Constraints (2h) imply that only the first individual can be replaced at time 1. Constraints (2i) impose the non-anticipativity constraint, forcing the decisions at the first stage to be the same. The constraints (2j) force all failed components at the first stage to be replaced. Constraints (2k) and (2l) define the auxiliary variable $Y_i^{r\omega}$. Constraints (2m) provide the definition of variable $y_{it}^{r\omega}$. Constraints (2n) – (2p) are the definition of variable $w_{it}^{r\omega}$. The remaining constraints (2q) – (2u) are binary constraints. The linearization of $\left|y_{it}^{r\omega}\right|$ can be found in Equation (4) and constraints (2v) – (2x).

### 4. Optimization algorithms

All decision variables in DEF are binary. Properties of stochastic integer programs are scarce, and general efficient methods are lacking. We therefore design a heuristic algorithm under the framework of PHA to solve practical-size problems within moderate CPU time. To assess the performance of the proposed heuristic algorithm in solving DEF, we compare the performance of the proposed algorithms with three conventional algorithms, namely, basic Benders decomposition (Algorithm 1), integer L-shaped method with Benders cuts (Algorithm 2) and standard PHA (Algorithm 3).

The basic Benders decomposition and integer L-shaped method with Benders cuts are considered because of LP relaxation and branch-and-cut are common methods for solving integer programs. Standard PHA (Watson and Woodruff 2011), which decomposes a problem by scenarios, provides a flexible framework for stochastic integer problem, and is also considered for comparison.

## 4.1 The Benchmark algorithms

### 4.1.1 Basic Benders decomposition algorithm

The basic Benders algorithm first solves the Benders master integer problem, and then solves the LP relaxation of subproblems to generate cuts which are added back to Benders master problem (Birge and Louveaux 2011). The procedure is repeated until no cuts found. We first define the initial master problem (MP) as follows:

$$\text{MP: min } dz + \sum_{i \in N} c_{i,\text{pr}} x_i + \sum_{i \in N} \left( c_{i,\text{cr}} - c_{i,\text{pr}} \right) \xi_i + \theta,$$

subject to:

$$\begin{aligned}
\theta &\geq Q(x), \\
x_i &\geq \xi_i, \quad i \in N \\
z &\geq x_i, \quad i \in N \\
x_i &\in \{0,1\}, \quad i \in N \\
z &\in \{0,1\}
\end{aligned}$$

where $Q(x) = \sum_{\omega \in \Omega} p(\omega) Q(x,\omega)$ and $Q(x, \omega)$ is the objective of the scenario $\omega$ in the second-stage problem, given by

$$Q(x,\omega) = \sum_{i \in N} \left( C_{i,\text{pr}} + C_{i,\text{cr}} - c_{i,\text{pr}} x_i - \left( c_{i,\text{cr}} - c_{i,\text{pr}} \right) \xi_i \right) + C_s - dz$$

and subject to constraints (2b) – (2x) except (2e), (2j) and (2r). Constraints (2e), (2j) and (2r) are excluded from the subproblem since they only concern the decision variables in the first stage. For each scenario $\omega$, Benders cut can be written as

$$\theta_\omega \geq \tilde{e}_{m,\omega} x - e_{m,\omega}, \quad m \in \{1, ..., M\}, \; \omega \in \Omega, \tag{5}$$

where $M$ denotes the maximum iterations.

---

**Algorithm 1: (Basic Benders decomposition)**

1: **Initialization:** $\theta_\omega \leftarrow -\infty$, for $\forall \omega \in \Omega$, $\tilde{\varepsilon} \leftarrow 10^{-2}$, and assign an integer feasible $x$ to the subproblem.

2: Solve the LP relaxation of the subproblem, $Q(x, \omega)$, for each $\omega \in \Omega$.

3: **If** $\theta_\omega - Q(x, \omega) \leq \tilde{\varepsilon} \; \forall \omega \in \Omega$, **return** optimal solution: $(x^*, \theta_\omega^*) \leftarrow (x, \theta_\omega)$. **Else**, go to step 4

4: Add Benders cuts using Equation (5) into the MP, where $\theta = \sum_{\omega \in \Omega} p(\omega) \theta_\omega$.

5: Solve the MP to get new $(x, \theta_\omega)$, $\forall \omega \in \Omega$. Go to step 2.

---

### 4.1.2. Integer L-shaped method with Benders cuts

In Algorithm 2, we initialize Benders master problem with Benders cuts. More specifically, the root node is obtained by solving the LP relaxation of the master problem via Benders decomposition and keeping the cuts. In the branch-and-cut process, at each node, if the solution is integer feasible, the subproblem is solved to generate integer optimality cuts which are defined as follows (Laporte and Louveaux 1993):

$$\theta \geq \left(Q(x^*) - L\right)\left(\sum_{i \in S(\mathbf{x}^*)} x_i - \sum_{i \notin S(\mathbf{x}^*)} x_i - |S(x^*)|\right) + Q(x^*), \tag{6}$$

where $S(x^*) := \{i \mid x_i^* = 1\}$.

In addition to the integer optimality cuts, Benders cuts are also generated and added into the MP if violated by the candidate solution, in order to improve the performance of the Integer L-shaped method. Therefore, for each node in the branch-and-cut search tree, if the candidate solution is integer feasible, both Benders cuts and integer optimality cuts are added via lazy constraint callback routine, otherwise only Benders cuts are added by using user-cut callback routine.

| **Algorithm 2: (Integer L-shaped method with Benders cuts)** |
|---|
| 1: **Initialization** $\theta^* \leftarrow +\infty$;<br>     Initialize the MP by solving the LP relaxation via Benders, and keep cuts $\Rightarrow (x, \theta)$<br><br>2: **Branch and Cut**<br>    At each node in the search tree:<br>        Solve LP relaxation $\Rightarrow (x, \theta)$<br>        **If** LP bound exceeds known incumbent $\theta^*$, prune.<br>        **If** x is integer feasible:<br>            Solve subproblem $Q(x)$ to generate integer optimality cuts using Equation (6).<br>            Solve LP relaxation of the subproblem $Q(x)$ to generate Benders cuts.<br>            **If** (**x**, $\theta$) violates any Benders cut or integer optimality cut, add cut to LP relaxation<br>                of the MP and resolve.<br>            **Else**, update the incumbent, $\theta^* \leftarrow \theta$<br>        **If** x is not integer feasible:<br>            Solve LP relaxation of the subproblem $Q(x)$ to generate Benders cuts.<br>            **If** (x, $\theta$) violates any Benders cut, add cut to LP relaxation and resolve.<br>            **Else**, branch to create new nodes. |

#### 4.1.3. Standard progressive hedging algorithm

We also examine the performance of the standard PHA on our problem. The PHA mitigates the computational difficulty associated with large problem instances by decomposing the extensive form according to scenario, and iteratively solving penalized versions of the subproblems to gradually enforce non-anticipativity (Aydin 2012). Solving individual scenario subproblems separately is generally much less computationally challenging and may allow a solver to exploit any special combinatorial structure that may be present. Moreover, the time expended for each iteration can be dramatically reduced by a very straightforward parallelization.

Specifically, the PHA proceeds by relaxing the non-anticipativity constraints using augmented Lagrangean relaxation and the problem becomes separable by each scenario. The scenario subproblems have augmented objective functions which include Lagrangean penalty functions corresponding to the relaxed non-anticipativity constraints. At each iteration of the PHA algorithm, these scenario subproblems are solved as deterministic problems. Solutions from all scenario sub-problems are then collected and averaged according to their non-anticipativity constraints and scenario probabilities. The deviation of each scenario sub-problem solution from these averages is used to update the Lagrangean multipliers. Next, the scenario sub-problems are re-solved with the updated augmented Lagrangean objective function. This iterative process continues until the Lagrangean dual problem converges to a solution satisfying the non-anticipativity constraints.

Details of the PHA are described in Algorithm 3. A different form of the objective function, $cx+E(Q(x, \omega))$, is used for a concise presentation of the algorithm (Gade et al. 2016).

| **Algorithm 3: (The standard PHA)** |
|---|
| 1. **Initialization**:<br>   Let $v \leftarrow 0$, $\tilde{\varepsilon} \leftarrow 10^{-2}$;<br>   $x_\omega^v \leftarrow \arg\min_\mathbf{x} (cx + Q(x,\omega)), \forall \omega \in \Omega$;<br>   $\bar{x}^v \leftarrow \sum_{\omega \in \Omega} p(\omega) x_\omega^v$;<br>   $w_\omega^v \leftarrow \rho(x_\omega^v - \bar{x}^v), \forall \omega \in \Omega$.<br>2. **Update the iteration counter**: $v \leftarrow v + 1$.<br>3. **Decomposition**:<br>   $x_\omega^v \leftarrow \arg\min_\mathbf{x}(cx + w_\omega^{v-1}x + \frac{\rho}{2}\|x - \bar{x}^{v-1}\| + Q(x,\omega)), \forall \omega \in \Omega$. |

4. **Aggregation**: $\bar{x}^v \leftarrow \sum_{\omega \in \Omega} p(\omega) x_\omega^v$.
5. **Update price**: $w_\omega^v \leftarrow w_\omega^{v-1} + \rho \left( x_\omega^v - \bar{x}^v \right)$, $\forall \omega \in \Omega$.
6. **Calculate converge distance**: $g^v \leftarrow \sum_{\omega \in \Omega} p(\omega) \left\| x_\omega^v - \bar{x}^v \right\|$, $\forall \omega \in \Omega$.
7. **Termination**: **If** $g^v < \tilde{\varepsilon}$, stop and **return** optimal solution $\bar{x}^v$. **Else**, go to step 2.

### 4.2. Progressive-hedging-based heuristic algorithm

Algorithm 1 cannot provide meaningful results due to the LP relaxation employed, and becomes more difficult and time-consuming as more cuts are added. Algorithm 2 is also computationally intensive as the number of binary variables and constraints increases. Standard PHA similarly suffers the computational intractability, since even for a small-scale multi-component maintenance problem, the scenario subproblem in the DEF can have a large number of decision variables and constraints, beyond what commercial solvers (e.g., CPLEX) can handle. However, PHA provides a flexible framework for solving stochastic integer problems. To address the bottleneck in solving the scenario subproblem using the standard PHA, i.e., step 3 in Algorithm 3, we develop an efficient heuristic algorithm based on the problem structure for the scenario subproblems, and use the PHA framework as a "wrapper" to force non-anticaptivity constraints.

The basic idea of our heuristic algorithm is as follows. Given a scenario subproblem, the heuristic iteratively groups maintenance activities of *working* individuals to reduce the setup costs and ultimately reduce the total maintenance costs. At each iteration it first obtains tentative replacement schedules, which are temporary, for all working individuals without considering economic dependence, and then considers a shifting window and groups maintenance activities within the shifting window. The tentative replacement schedule and the shift window are optimized to find the lowest total maintenance cost. The time complexity of the heuristic algorithm is polynomial (see Appendix A.1 for proof).

Before describing the details of the heuristic, we first present two properties regarding the optimal solution.

**Theorem 1.** For each scenario subproblem, there exists an optimal solution such that at each decision period $t \in T_s$, if there is any group of individuals (including one-individual group) that is maintained, there is at least one individual in the group that is replaced at one time unit before its failure or at its failure, except for all types of components' last individuals that are replaced in the planning horizon. (Proof is in Appendix A.2). □

**Theorem 2.** Given a set of working individuals sorted according to their failure times, there exists an optimal solution for this set such that maintenance activities are executed following the same order. (Proof is shown in Appendix A.3). □

Theorem 1 helps determine tentative replacement schedules for each individual. Based on Theorem 1, we only need to consider two tentative replacement schedules for individuals when ignoring economic dependence, i.e., replacing onetime unit before a failure or at the failure. Theorem 2 further ensures that it is optimal to execute the replacement activities for all working individuals in the order they are tentatively planned according to Theorem 1. Theorems 1 and 2 significantly decrease the number of possible feasible solutions needed to be considered in the heuristic, and thus substantially reduce the algorithm complexity.

The details of the heuristic are described as follows. Let $K$ denote the set of all working individuals at the current iteration of the heuristic. Let $\beta_{ij}$ and $\beta'_{ij}$ denote the tentative and actual replacement times of individual $I_{ij}$, respectively. Let $K'$ represent the sorted set of $K$ according to $β_{ij}$, and $K'[i]$ represent the individual in the $i^{th}$ position in set $K'$. For example, consider a four-component system, and the tentative replacement times of four individuals at one iteration are provided in Figure 4. In this example, we have $K = \{I_{1,5}, I_{2,3}, I_{3,2}, I_{4,4}\}$ and $K' = \{I_{2,3}, I_{1,5}, I_{4,4}, I_{3,2}\}$.

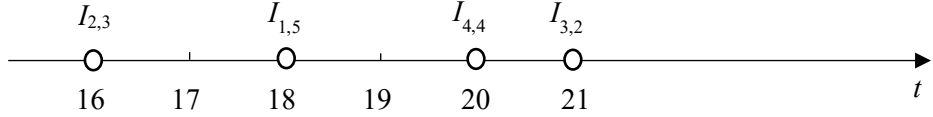

Figure 4: Working individuals at one iteration

Given a shifting window $\iota$ and a working individual set $K$, the **Grouping Rule** constructs candidate grouping options and selects the one with the minimum weighted PHA replacement cost cumulated till the current iteration. Specifically, the first candidate grouping option starts from $K'[1]$ and groups all current working individuals with a tentative replacement time between $β_{K'[1]}$ and $β_{K'[1]} + \iota$. Let $K'[v]$ represent the last individual grouped with $K'[1]$. The next shifting window starts from $K'[v+1]$ and ends at $K'[v+1]+\iota$. Working individuals with a tentative replacement time in this window are grouped together. The grouping process of the first option continues until $v = |K'|$, implying no more individual can be grouped. The second grouping option starts from $K'[2]$ and the same grouping process is repeated. The total number of group options is $|K'| - 1$. Note that in any grouping option, if individual $K'[1]$ is not grouped with other individual(s), it will become a one-individual group and replaced by its immediate successor in order to keep the grouping process rolling in the time horizon. The tentative replacement time of any individual that is grouped with one or more individuals in the optimal grouping option for set $K$ becomes its actual replacement time and is replaced with its immediate successor. Note that if the new individual has a tentative replacement schedule beyond the planning horizon, it is removed from set *K*. Now we have a new set of working individuals and the **Grouping Rule** is applied for the new set. The heuristic stops when set *K* is empty.

We use the same four-component system considered earlier to illustrate the grouping process using the **Grouping Rule**. Suppose the shifting window $\iota = 3$. The three candidate group options are illustrated in Figure 5. Among all three options, group option 3 has the minimum weighted PHA cost cumulated. Therefore, the actual replacement time for $I_{2,3}$ is 16 and 20 for $I_{4,4}$ and $I_{3,2}$.

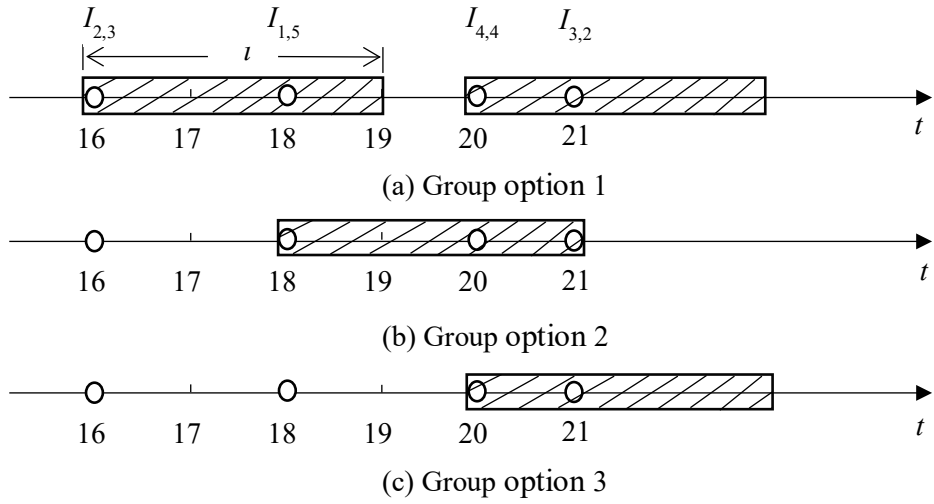

Figure 5: Group options at one iteration

The heuristic is summarized in Algorithm 4 and the Grouping Rule.

---

**Algorithm 4: (Heuristic algorithm for one scenario)**

**Initialization**: not-used-residual lifetime $\Delta \leftarrow \{0, 1\}$, and determine a set of values for $\iota$, $\iota = \{\iota_1, \iota_2, \ldots\}$

For all combinations of $\Delta$ and $\iota$, select the one that has the minimum total weighted PHA replacement cost and return the corresponding optimal solution

1: **Initialization**: Assign tentative replacement times for the first individual of each component
$$K \leftarrow \{I_{1,1}, I_{2,1}, \ldots, I_{n,1}\}, \text{ and } \beta_{i,1} \leftarrow T_{i,1} - \Delta, \forall i \in N.$$
2: Apply **Grouping Rule** to obtain the optimal group option $W$.
3: Update Set $K$.
$\forall I_{ij} \in W$
  Replace $I_{ij}$ in set $K$ with $I_{i,j+1}$;
  Assign tentative replacement schedule to $I_{i,j+1}$, $\beta_{i,j+1} \leftarrow \beta'_{ij} + (T_{i,j+1} - \Delta)$;
  $\forall I_{ij} \in K$, **If** $\beta_{ij} > T$, remove $I_{ij}$ from set $K$;

| **If** $K$ is empty, stop. **Else**, go to step 2. |
|---|

| **Grouping Rule** |
|---|
| 1: Sort $K$ in ascending order based on $\beta_{ij}$ ⇒ sorted set $K'$. |
| 2: Select the option $m$ that has the lowest weighted PHA replacement cost cumulated ⇒ $W$.<br>   **Group option $m$**: $m$ from 1 to $|K'| - 1$<br>        2.1: Initialize the last individual grouped as $K'[\upsilon]$: $\upsilon \leftarrow m$<br>        2.2: $t \leftarrow \beta(K'[\upsilon])$<br>            Group individuals in set $K'$ if the actual replacement times of their predecessors are before $t$ until $\beta(K'[\upsilon']) > t + \iota$, $\upsilon' = \upsilon + 1, \upsilon + 2, \ldots$<br>        2.3: Let $\vartheta$ denote the position of the last individual grouped in step 2.2<br>            Update actual replacement times: $\tau'(K'[\upsilon']) \leftarrow t$, $\upsilon' = \upsilon + 1, \upsilon + 2, \ldots, \vartheta$<br>            $\upsilon \leftarrow \vartheta + 1$<br>        2.4: **If** $\upsilon \geq |K'|$, compute the weighted PHA replacement cost cumulated, then stop.<br>           **Else**, go to step 2.2. |
| 3: **If** $K'[1] \notin W$: $W \leftarrow W \cup \{(K'[1])\}$. |
| 4: **Return** set $W$. |

Our modeling and solution approach can be easily extended to incorporate several other modeling aspects, such as cost discounts depending on the number of components maintained, restriction on the number of components replaced together, and nonlinear operation costs. It can also be easily adapted to solve opportunistic maintenance (OM). We can find the optimal OM policy that groups maintenance only at CM by restricting the not-used-residual lifetime in our heuristic to $\Delta = 0$, and find the optimal OM policy that only groups maintenance at PM by only allowing $\Delta = 1$. More structured policies can be derived based on the patterns of results from our general model. For example, fixed groups in direct-grouping can be obtained based on how frequent some components are grouped for maintenance in our results. Indirect grouping policy can be derived based on the individual PM replacement intervals.

## 5. Computational study

In this computational study, we first introduce an analytical approach to determining the minimum number of scenarios needed for a given accuracy, then examine the performance of all four algorithms, and assess the approximation performance of the multi-stage model by using proposed heuristic algorithm (Algorithm 4) in a rolling horizon. We perform our computational study on a computer with a CPU of Intel i7-6700, 3.4G Hz and a RAM of 16G. A python based package Pyomo (Hart et al. 2011, Hart et al. 2012) is used to implement the algorithms with the solver of CPLEX v12.7.1.

### 5.1 Scenario sampling

Consider a stochastic programming problem $V^* = \min_{x \in X}\{f(x) := E[F(x, \xi)]\}$ and its sample average approximation (SAA) problem $\hat{V}_{|\Omega|} = \min_{x \in X}\{\hat{f}_{|\Omega|}(x) := \frac{1}{|\Omega|}\sum_{i=1}^{|\Omega|} F(x, \xi^i)\}$, where $X$ is the feasible region of decision variable $x$, and $\Omega$ is a sample of the random vector $\xi$. Denote the $\varepsilon$-optimal solution sets of the stochastic program and SAA problem by $S^\varepsilon := \{x \in X : f(x) \leq V^* + \varepsilon\}$ and $\hat{S}^\varepsilon_{|\Omega|} := \{x \in X : \hat{f}_{|\Omega|}(x) \leq \hat{V}_{|\Omega|} + \varepsilon\}$ respectively. To guarantee that the optimal solution of SAA problem is the $\varepsilon$-optimal solution of the true stochastic programming problem with probability $1 - \alpha$, Shapiro et al. (2009) provide an analytical expression of the minimum sample size required. Theorem 3 summarizes the main result in their study.

**Theorem 3 (Theorem 5.17 (Shapiro et al. 2009)).** Suppose there exists a constant $\sigma > 0$ such that for any $x \in X \backslash S^\varepsilon$ the moment generating function $M_x(t)$ of the random variable $H(x, \xi) - E[H(x, \xi)]$ satisfies $M_x(t) \leq e^{(\sigma^2 t^2/2)}, \forall t \in \mathbb{R}$, then for $\varepsilon > 0$, $0 \leq \tau < \varepsilon$ and $\alpha \in (0,1)$, and sample size $|\Omega|$ satisfying

$$|\Omega| \geq \frac{2\sigma^2}{(\varepsilon-\tau)^2}\ln(\frac{|X|}{\alpha}), \quad (7)$$

it follows that

$$\Pr(\hat{S}^\tau_{|\Omega|} \subset S^\varepsilon) \geq 1-\alpha,$$

where $H(x,\xi) := F(u(x),\xi) - F(x,\xi)$, $x \in X \setminus S^\varepsilon$ and mapping $u : x \in X \setminus S^\varepsilon \to X$ satisfies $f(u(x)) \leq f(x) - \varepsilon^*$ for some $\varepsilon^* \geq \varepsilon$. □

To determine constant $\sigma$, Shapiro et al. (2009) show that if $|H(x,\xi) - E[H(x,\xi)]| < b$ is satisfied with some $b > 0$ for all $x \in X$, then $\sigma^2 := b^2$. If $\varepsilon^*$ is small compared to $\max_x F(x,\xi)$, then in DEF, an upper bound $b$ can be estimated by considering CR for all individuals over the planning horizon.

$$|H(x,\xi) - E[H(x,\xi)]| \leq |H(x,\xi)| + |E[H(x,\xi)]| \approx |H(x,\xi)|$$

$$\leq |F(u(x),\xi)| + |F(x,\xi)| \leq 2F(x,\xi) \leq 2T(\sum_{i \in N} c_{i,\text{cr}} + d) = b \quad (8)$$

### 5.2 Performance of algorithms in solving DEF

In this section, we compare the computational times and cost errors of the four algorithms in solving DEF. For standard PHA (Algorithm 3), we run our experiments in a stochastic programming package PySP inside Pyomo. Assume that all components' lifetimes follow Weibull distributions. For each component, we draw the shape and scale parameters from uniform distributions $U(4, 7)$ and $U(1, 8)$, respectively. The cost of CR ($c_{i,\text{cr}}$) is drawn from a uniform distribution $U(6, 16)$. Parameter values that are randomly drawn and used in the computational study is provided in Appendix A.4. Without loss of generality, the cost of PR ($c_{i,\text{pr}}$) is assumed to be 1. Suppose that setup cost $d$ is 5, the initial ages of all first individuals are 2, and the individual of component 1 is assumed to be failed at the first stage.

For a fixed planning horizon $T = 10$, Table 1 shows the problem size of each scenario subproblem with different number of components. The number of constraints is instance-dependent because of constraints (2f) and (2m), and we approximate it with the maximum possible number of constraints. From Table 1, we can see that the problem size is very large even for a small problem.

Table 1: Illustration of problem size

| $n$ | Variables | | Constraints | |
|---|---|---|---|---|
| | master problem | subproblem | master problem | subproblem |
| 2 | 3 | 3,970 | 6 | 10,575 |
| 4 | 5 | 7,930 | 12 | 21,139 |
| 6 | 7 | 11,890 | 18 | 31,703 |
| 8 | 9 | 15,850 | 24 | 42,267 |

The number of scenarios needed for each test case in Table 3 is summarized in Table 2, which is determined using Equations (7) and (8) by choosing $\varepsilon = 0.1\sigma$, $\tau = 0.1\varepsilon$ and $\alpha = 0.1$, This parameter setting guarantees the optimal solution of SAA problem is an $0.1\sigma$-optimal solution of the true stochastic programming problem with probability 0.9. We use the same parameter setting to determine the number of scenarios throughout the paper.

Table 2: Number of scenarios needed

| $n$ | $|\Omega|$ | $n$ | $|\Omega|$ |
|---|---|---|---|
| 2 | 740 | 5 | 1250 |
| 3 | 910 | 6 | 1420 |
| 4 | 1080 | 7 | 1600 |

We compare the performance of Algorithms 1 – 4 for 18 cases. Table 3 summarizes the performance of the four algorithms. NA is reported if the computational time is longer than one day or out of memory, or if the true objective of DEF is not available for computing objective percentage error. From Table 3, we can see that Algorithms 1 – 3 can only solve small problems (e.g., $n \leq 4$) and Algorithm 4 is the only algorithm that can solve all test cases efficiently. We further examine the performance of Algorithm 4. We compute the percentage error between the objective values from using Algorithm 4 and solving DEF exactly by CPLEX. We can see that the performance of Algorithm 4 is acceptable for small problems with a maximum percentage error of 12.73%. Based on our computational studies, the proposed heuristic algorithm (Algorithm 4) performs well comparing to the true objective value of the DEF model and is capable of solving practically large-scale problems.

Table 3: Algorithm performance in solving DEF

| | | | solver | | Algorithm 1 | | | | Algorithm 2 | | | Algorithm 3 | | | Algorithm 4 | | | |
|---|---|---|---|---|---|---|---|---|---|---|---|---|---|---|---|---|---|---|
| case | $n$ | $T$ | CPU time (sec.) | Obj. | Iterations | Cuts | CPU time (sec.) | Obj. | Cuts | CPU time (sec.) | Obj. | Iterations | CPU time (sec.) | Obj. | Iterations | CPU time (sec.) | Obj. | Obj. error % |
| 1 | 2 | 6 | 103 | 25.80 | 3 | 1478 | 466 | -99.65 | 4 | 1156 | 25.80 | 3 | 1361 | 25.80 | 1 | 17 | 26.67 | 3.37% |
| 2 | 2 | 8 | 202 | 28.63 | 3 | 1474 | 546 | -111.53 | 4 | 1556 | 28.63 | 2 | 1315 | 28.63 | 1 | 14 | 30.92 | 8.00% |
| 3 | 2 | 10 | 503 | 32.82 | 3 | 1380 | 760 | -167.22 | 4 | 1670 | 32.82 | 3 | 1775 | 32.82 | 2 | 25 | 34.88 | 6.28% |
| 4 | 3 | 6 | 237 | 26.92 | 3 | 1818 | 884 | -143.87 | 9 | 2505 | 26.92 | 5 | 2502 | 26.92 | 2 | 51 | 28.55 | 6.45% |
| 5 | 3 | 8 | 561 | 30.14 | 4 | 1830 | 1613 | -187.43 | 13 | 4070 | 30.14 | | NA | | 6 | 103 | 33.31 | 10.52% |
| 6 | 3 | 10 | 2777 | 34.50 | 5 | 3070 | 2575 | -239.20 | | NA | | | NA | | 4 | 92 | 37.94 | 9.97% |
| 7 | 4 | 6 | 472 | 30.02 | 5 | 3263 | 2198 | -176.37 | | | | | | | 4 | 145 | 33.84 | 12.72% |
| 8 | 4 | 8 | NA | | 6 | 4764 | 3895 | -227.89 | | NA | | | NA | | 4 | 181 | 40.02 | NA |
| 9 | 4 | 10 | NA | | | | NA | | | | | | | | 4 | 269 | 46.05 | NA |
| 10 | 5 | 6 | | | 5 | 2848 | 3052 | -221.73 | | | | | | | 4 | 324 | 39.97 | |
| 11 | 5 | 8 | NA | | | | NA | | | NA | | | NA | | 4 | 422 | 48.30 | NA |
| 12 | 5 | 10 | | | | | NA | | | | | | | | 4 | 551 | 56.66 | |
| 13 | 6 | 6 | | | | | | | | | | | | | 4 | 623 | 44.9 | |
| 14 | 6 | 8 | NA | | | | NA | | | NA | | | NA | | 4 | 837 | 55.16 | NA |
| 15 | 6 | 10 | | | | | | | | | | | | | 4 | 1117 | 65.35 | |
| 16 | 7 | 6 | | | | | | | | | | | | | 4 | 1385 | 53.64 | |
| 17 | 7 | 8 | NA | | | | NA | | | NA | | | NA | | 5 | 2812 | 67.09 | NA |
| 18 | 7 | 10 | | | | | | | | | | | | | 5 | 3398 | 80.47 | |

**5.3 Performance in approximating the multi-stage model using a rolling horizon approach**

In this section, we compare our heuristic in a rolling horizon with a direct-grouping approach (Van Horenbeek and Pintelon (2013)). We also compare the performances using our approach in a rolling horizon with the results from solving the multi-stage model using an exact method, and assess the performance of the proposed approximation method.

The direct-grouping model (Van Horenbeek and Pintelon 2013) uses a dynamic-programming algorithm that first finds the optimal replacement schedule for each component without considering economic dependence and then sort the components based on that. At iteration $j$, the algorithm identifies two groups that cover all maintenance activities of components 1 to $j$ and provide the best savings for these components. The best grouping structure can be found by backtracking. This algorithm has in the worst case a time complexity of $o(n^2)$. However, the limitation of this algorithm is that it only considers the group structure of two groups at each iteration and ignores all other options (e.g., partition all maintenance activities into three or more groups).

We conduct a sensitivity analysis to assess the performance of Algorithm 4. Suppose the length of a decision period $\delta$ equals to 1. At each decision period, we consider a two-stage stochastic maintenance optimization problem (DEF) where the second stage combines decisions of the remaining periods. We repeat this procedure 5 times to obtain the average total maintenance cost over the planning horizon. The PR cost is assumed to be 1 for each type of component, and the CR costs are drawn from two different uniform distributions, $U(6, 16)$ and $U(17, 27)$. The lifetime of each individual is assumed to follow a Weibull distribution. To introduce more heterogeneity to the system, the shape parameter of the Weibull distribution is drawn from two uniform distributions $U(1, 3)$ and $U(4, 7)$, and the scale parameter of the Weibull distributions is drawn from two different distributions, $U(1, 5)$ and $U(5, 10)$. Two levels of setup costs are considered. Assume all working individuals are functioning and have an age of 0 at the first decision period. Different levels of parameters are provided in Table 4. Parameter values that are randomly drawn and used for the comparison are summarized in Appendix A.5.

Table 4: Different levels of parameter

| Level | shape parameter | scale parameter | $d$ | cost of CR |
|---|---|---|---|---|
| High | $U(4,7)$ | $U(5, 10)$ | 100 | $U(17, 27)$ |
| Low | $U(1,3)$ | $U(1, 5)$ | 5 | $U(6, 16)$ |

Table 5 summarizes the comparison results with different number of components $n$ and planning horizon $T$. From Table 5, we can see that solving our two-stage model by the proposed heuristic in a rolling horizon provides good approximation of the multi-stage problem for all cases examined and the heuristic significantly outperforms the benchmark algorithm. In particular, our approach shows a better performance when the setup cost $d$ is high.

Table 5: Numerical example for rolling horizon comparison ($n = 2$, $T = 10$, $|\Omega| = 910$)

| case | shape | scale | $d$ | CR cost | multi-stage (P1) | Benchmark obj. | error % | Algorithm 4 obj. | error % |
|---|---|---|---|---|---|---|---|---|---|
| 1 | H | H | H | H | 105.39 | 212.76 | 101.90% | 123.99 | 17.70% |
| 2 | H | H | H | L | 104.51 | 208.36 | 99.40% | 112.99 | 8.10% |
| 3 | H | H | L | H | 8.88 | 27.64 | 211.30% | 11.18 | 25.90% |
| 4 | H | H | L | L | 8 | 21.04 | 163.00% | 9.05 | 13.10% |
| 5 | H | L | H | H | 268.3 | 434.88 | 62.10% | 281.84 | 5.00% |
| 6 | H | L | H | L | 244.6 | 421.68 | 72.40% | 276.12 | 12.90% |
| 7 | H | L | L | H | 31.53 | 54.88 | 74.10% | 35.43 | 12.40% |
| 8 | H | L | L | L | 28.03 | 45.84 | 63.50% | 33.44 | 19.30% |
| 9 | L | H | H | H | 136.13 | 417.76 | 206.90% | 149.08 | 9.50% |
| 10 | L | H | H | L | 123.43 | 310.76 | 151.80% | 135.48 | 9.80% |
| 11 | L | H | L | H | 22.48 | 46.32 | 106.00% | 24.12 | 7.30% |
| 12 | L | H | L | L | 16.79 | 31.12 | 85.30% | 18.6 | 10.80% |
| 13 | L | L | H | H | 340.36 | 958.92 | 181.70% | 352.2 | 3.50% |
| 14 | L | L | H | L | 306.54 | 938.12 | 206.00% | 319 | 4.10% |
| 15 | L | L | L | H | 60.12 | 125.92 | 109.40% | 62.24 | 3.50% |

| 16 | L | L | L | L | 44.7 | 83.12 | 86.00% | 49.16 | 10.00% |

We further compare the performance of Algorithm 4 and the benchmark algorithm on large problems. The results are summarized in Table 6. Note that exact method cannot solve the problem instances in Table 6 and are therefore not included. Similar to what we observed in Table 5, when the setup cost is high, Algorithm 4 leads to more cost savings.

Table 6: Numerical example for rolling horizon comparison (part 1, $n = 4$, $T = 10$, $|\Omega| = 1{,}250$)

| case | shape | scale | $d$ | CR cost | Benchmark obj. | Algorithm 4 obj. | savings |
|---|---|---|---|---|---|---|---|
| 1 | H | H | H | H | 216.1 | 119.05 | 97.05 |
| 2 | H | H | H | L | 210.3 | 110.25 | 100.05 |
| 3 | H | H | L | H | 34.2 | 23.96 | 10.24 |
| 4 | H | H | L | L | 24.6 | 15.25 | 9.35 |
| 5 | H | L | H | H | 464.7 | 363.44 | 101.26 |
| 6 | H | L | H | L | 429.92 | 325.64 | 104.28 |
| 7 | H | L | L | H | 76.82 | 52.44 | 24.38 |
| 8 | H | L | L | L | 49.92 | 43.64 | 6.28 |
| 9 | L | H | H | H | 433.18 | 225.24 | 207.94 |
| 10 | L | H | H | L | 320.86 | 203.24 | 117.62 |
| 11 | L | H | L | H | 61.06 | 58.84 | 2.22 |
| 12 | L | H | L | L | 39.56 | 36.46 | 3.1 |
| 13 | L | L | H | H | 1021.22 | 480.28 | 540.94 |
| 14 | L | L | H | L | 972.02 | 483.58 | 488.44 |
| 15 | L | L | L | H | 166.38 | 103.66 | 62.72 |
| 16 | L | L | L | L | 117.02 | 83.88 | 33.14 |

Table 6 (cont'd): Numerical example for rolling horizon comparison (part 2, $n = 8$, $T = 20$, $|\Omega| = 1{,}940$)

| case | shape | scale | $d$ | CR cost | Benchmark obj. | Algorithm 4 obj. | savings |
|---|---|---|---|---|---|---|---|
| 17 | H | H | H | H | 557.5 | 368.56 | 188.94 |
| 18 | H | H | H | L | 543.8 | 352.52 | 191.28 |
| 19 | H | H | L | H | 103.82 | 94.02 | 9.8 |
| 20 | H | H | L | L | 78.92 | 68.26 | 10.66 |
| 21 | H | L | H | H | 2155.04 | 1173.36 | 981.68 |
| 22 | H | L | H | L | 2090.46 | 986.74 | 1103.72 |
| 23 | H | L | L | H | 393.2 | 236.34 | 156.86 |
| 24 | H | L | L | L | 284.36 | 212.96 | 71.4 |
| 25 | L | H | H | H | 1029.48 | 887.54 | 141.94 |
| 26 | L | H | H | L | 712.86 | 680.6 | 32.26 |
| 27 | L | H | L | H | 297.32 | 211.98 | 85.34 |
| 28 | L | H | L | L | 153.64 | 137.36 | 16.28 |
| 29 | L | L | H | H | 2445.58 | 2076.54 | 369.04 |
| 30 | L | L | H | L | 2226.4 | 1535.12 | 691.28 |
| 31 | L | L | L | H | 640.58 | 563.5 | 77.08 |
| 32 | L | L | L | L | 421.4 | 379.3 | 42.1 |

## 6. Conclusion and future research

In this paper, we consider the problem of multi-component maintenance optimization over the finite planning horizon. We formulate the problem as a multi-stage decision-dependent stochastic integer program, and approximate it with a novel two-stage stochastic linear integer model in a rolling horizon. The proposed models are general with no restrictions on maintenance grouping. A progressive-hedging-based heuristic is designed to solve practically large-size two-stage problems. To assess the performance of the heuristic, we compare it with three conventional algorithms and our computational studies show that the proposed heuristic provides satisfying results and is capable of solving practically large-scale problems. We also evaluate the performance of the heuristic in a rolling horizon relative to the true global optimal for small problems. Results show that solving our two-stage model by the proposed heuristic in a rolling horizon provides a good approximation of the multi-stage problem. The proposed heuristic in a rolling horizon is further benchmarked with a widely studied dynamic-programming-based algorithm. Our heuristic significantly outperforms the benchmark algorithm on all cases examined.

Our work has extended the available literature in multi-component maintenance by using stochastic programming approach. The modeling and solution techniques developed in this paper opens new research and implementation opportunities. Future research will consider a different widely used maintenance policy, condition-based maintenance (CBM). CBM leverages sensor information on components' health status through inspection or real-time monitoring and aims to perform maintenance just in time by setting optimal control thresholds. Capturing these complexities requires a different problem formulation and different optimization algorithms. Moreover, maintenance activities are often subject to a pre-determined budget with a requirement on a system's reliability or availability. Future work will incorporate these constraints into the decision model. Lastly, it is worth extending the problem for more complex systems with stochastic and structural dependences, in addition to the economic dependence.

**Appendix**

**(A.1) The time complexity of the heuristic algorithm (Algorithm 4) is polynomial.** □

**Proof:**

The inputs that related to the time complexity are (1) not-used-residual lifetimes of an individual $|\Delta|$, (2) shifting windows $|\iota|$, (3) length of planning horizon $T$ and (4) number of component $n$.

Denote $G(\cdot)$ as the running time function.

The running time regarding different $n$ can be evaluated by the number of group options found in Grouping Rule. For the 1st group option, it requires $(n-1)$ steps. For the 2nd group option, it requires $(n-2)$ steps. So forth, for the $(n-1)$ group option that found with a given $n$, the total steps are $(n-1) + (n-2) + \ldots + 1 = n(n-1)/2$. Therefore, the running time in terms of $n$ is

$$G(n) = n(n-1)/2 + (n-1)c_1$$

where $c_1$ is a constant that represents the execution time of other statement in Grouping Rule.

In the worst case, the Grouping Rule is executed at every time point, i.e. $T$ times. Thus, the running time in terms of $n$ and $T$ is

$$G(n,T) = TG(n) + Tc_2$$

where $c_2$ is a constant that represents the execution time of other statement.

Because $\Delta$ and $\iota$ are the search variables, the total running time is

$$G(n,T,|\Delta|,|\iota|) = |\Delta| \cdot |\iota| \cdot G(n,T) = n^2 T |\iota| + nT|\iota|c_3 + T|\iota|c_4$$

where $c_3 = 2c_1 - 1$ and $c_4 = 2(c_2 - c_1)$. Notice that $|\Delta| = 2$ based on Theorem 1. Using big-$O$ notation, the time complexity of this algorithm is $O(n^2 T |\iota|)$, which is polynomial.

Proof completed. □

**(A.2) Theorem 1.** For each scenario subproblem, there exists an optimal solution such that at each decision period $t \in T_s$, if there is any group of individuals (including one-individual group) that is maintained, there is at least one individual in the group that is replaced at one time unit before its failure or at its failure, except for all types of components' last individuals that are replaced in the planning horizon. □

**Proof:**

Let group $m$ be the set of individuals that replaced at time $t_m$, and group structure $W_t$ collect all groups that replaced at or before time $t \leq T$ in chronological order. Let $\varepsilon_m$ denote the minimum not-used-residual lifetime of all individuals in group $m$ in a group structure $W_T$ in the planning horizon $T_s$. Let $W_{last}$ include the groups such that in each of these groups, every individual is the last individual of its component type that replaced in the planning horizon $T_s$. We need to prove that there exists an optimal grouping structure $W_T^*$ such that $\varepsilon_m \leq 1, \forall m \in W_T^* \setminus W_{last}$.

We prove the theorem by contradiction. Suppose that there exists at least one group satisfying $\varepsilon_m > 1$ in any optimal group structure $W_T^*$. The goal is to show that by appropriately regrouping the maintenance activities, we can find an alternative group structure $W_T'$ that yields no higher cost and has $\varepsilon_m \leq 1$, $\forall m \in W_T' \setminus W_{last}$

We start the proof by first considering the case where there is only one group, say group $\lambda$, with $\varepsilon_\lambda > 1$ and $\varepsilon_m \leq 1$, $\forall m \in W_T' - W_{last} - \{\lambda\}$. If there are more than one groups of this kind, we will start with the last group with $\varepsilon_m > 1$ and perform the regrouping process iteratively until $\varepsilon_m \leq 1$ for all groups.

Suppose groups $\lambda \notin W_{last}$ and $(\lambda+1) \notin W_{last}$ are replaced at $t_\lambda$ and $t_{\lambda+1}$ in group structure $W_T^*$ respectively. Next, we describe the details of how we construct the alternative group structure $W_T'$. We use the auxiliary variable $w_{it}^r$ to help building the new group structure. We construct the new group structure by shifting and regrouping individuals in group $\lambda$ and the subsequent groups when needed. There are three possible scenarios in the regrouping process.

**Scenario 1:** $t_{\lambda+1} \geq t_\lambda + \varepsilon_\lambda$.

In this scenario, we let all individuals in group $\lambda$ replaced at $t_\lambda + \varepsilon_\lambda - 1$. We have $\left(w_{it}^r\right)'$ given by

$$\left(w_{it}^r\right)' = \begin{cases} w_{it}^r, & I_{ir} \notin \lambda, t \in T_s \\ 1, & I_{ir} \in \lambda, t = t_\lambda + \varepsilon_\lambda - 1 \\ 0, & I_{ir} \in \lambda, t \neq t_\lambda + \varepsilon_\lambda - 1 \end{cases} \quad \text{(A1)}$$

Based on $\tilde{x}_{it'}^r = \sum_{t=1}^{t'} w_{it}^r$, we have,

$$\left(\tilde{x}_{it}^{r}\right)' = \begin{cases} \tilde{x}_{it}^{r}, & I_{ir} \notin \lambda, t \in T_s \\ 0, & I_{ir} \in \lambda, t < t_\lambda + \varepsilon_\lambda - 1 \\ 1, & I_{ir} \in \lambda, t \geq t_\lambda + \varepsilon_\lambda - 1 \end{cases}. \tag{A2}$$

From the definition of $z_t$, we have

$$\left(z_t\right)' = \begin{cases} z_t, & t \neq t_\lambda, t \neq t_\lambda + \varepsilon_\lambda - 1 \\ 0, & t = t_\lambda \\ 1, & t = t_\lambda + \varepsilon_\lambda - 1 \end{cases}. \tag{A3}$$

With the Equations (A1) – (A3), it is straightforward that solution $\left(w_{it}^{r}\right)'$ does not violate any constraint. It can also be easily verified that the cost remains the same in the reconstructed group structure under this scenario. The new replacement schedule of group $\lambda$ make $\varepsilon_\lambda \leq 1$, but it may cause some subsequent group(s) to have the minimum not-used-residual lifetime of all individuals in the group greater than one. We will perform this regrouping process recursively until all groups satisfy $\varepsilon_m \leq 1$.

**Scenario 2:** $t_{\lambda+1} < t_\lambda + \varepsilon_\lambda$.

Let $\varepsilon_{ir}$ denote the not-used-residual life of individual $I_{ir}$. Define set $P$ such that $\forall I_{ir} \in P$, we have $I_{ir} \in \lambda+1, I_{i,r-1} \notin \lambda$ and $\varepsilon_{ir} \leq 1$. We further separate Scenario 2 into two sub-scenarios based on whether set $P$ is empty.

**Scenario 2.1:** $P \neq \varnothing$.

Define set $Q$ such that $\forall I_{ir} \in Q$, we have $I_{ir} \in \lambda+1, I_{i,r-1} \in \lambda$. We construct two new groups $\lambda'$ and $(\lambda+1)'$, such that $\lambda' = \lambda \cup (\lambda+1) - Q$ and $(\lambda+1)' = Q$.

We replace all individuals in group $\lambda'$ at time $t_{\lambda+1}$. It is obvious that the new group $\lambda'$ satisfies $\varepsilon_{\lambda'} \leq 1$. We then process group $(\lambda+1)'$ in the same way as how we process group $\lambda$ from the beginning. It is also obvious that no additional cost incurs during this regrouping. And we will regroup recursively until all groups satisfy $\varepsilon_m \leq 1$.

If set $Q$ is empty, we actually combined two groups into one group and satisfy $\varepsilon_m \leq 1$ for all groups while saving one setup cost.

We can similarly show that no constraint is violated during the regrouping process in this scenario, and detailed proof is omitted.

**Scenario 2.2:** $P = \emptyset$

We construct the same two new groups $\lambda'$ and $(\lambda + 1)'$ as in Scenario 2.1. The difference is that the new group $\lambda'$ in this scenario does not satisfy $\varepsilon_{\lambda'} \leq 1$ at replacement time $t_{\lambda+1}$. We next process group $\lambda'$ in the same way as how we process group $\lambda$ from the beginning and then do the same thing for group $(\lambda + 1)'$.

Proof completed. □

**(A.3) Theorem 2.** Given a set of working individuals sorted according to their failure times, there exists an optimal solution for this set such that maintenance activities are executed following the same order. □

**Proof:**

We prove Theorem 2 by contraction.

Let group $m$ be the set of individuals that replaced at time $t_m$, and group structure $W_t$ collect all groups that replaced at or before time $t \leq T$ in chronological order. For two individuals $I_{ir}$ and $I_{jr}$ ($i \neq j$) in a working individual set, denote their failure times as $t_1$ and $t_2$ ($t_1 < t_2$) respectively. Let $t'_1$ and $t'_2$ represent the actual replacement times for individuals $I_{ir}$ and $I_{jr}$, respectively. Suppose $t'_1 = \eta_1$ and $t'_2 = \eta_2$, and there is an optimal structure $W^*$ that has $\eta_1 > \eta_2$. The goal is to show that we can find a new optimal group structure which satisfies $t'_1 \leq t'_2$ at the same or a lower cost.

Consider a group structure $W'$ where both individuals are replaced at $t'_1$. It is obvious that individual $I_{jr}$ is preventively replaced in $W^*$, and delay its replacement time to $\eta_1$ does not change its replacement type, meaning no additional cost because of no change in the replacement type. If individual $I_{jr}$ is grouped with some other individuals in $W^*$, it is obvious that the cost of group structure $W'$ is the same as $W^*$. If individual $I_{jr}$ is not grouped with any other individual in $W^*$, then we eliminate one setup cost in group structure $W'$, which leads to a lower cost.

Proof completed. □

**(A.4)**

Table A.4. Component parameters used for Table 3

| $i$ | shape | scale | $c_{i,cr}$ |
|---|---|---|---|
| 1 | 6.5 | 6.9 | 14.4 |
| 2 | 6.7 | 5 | 11.4 |
| 3 | 5.4 | 7.3 | 9.4 |
| 4 | 4.9 | 4.8 | 8.0 |
| 5 | 4.8 | 4.2 | 11.1 |
| 6 | 4.4 | 4.5 | 14.2 |
| 7 | 5.5 | 3.2 | 7.4 |

**(A.5)**

Table A.5: Parameters for each component in different levels used for Table 5

| $i$ | shape parameter | | scale parameter | | $c_{i,cr}$ | |
|---|---|---|---|---|---|---|
| | High | Low | High | Low | High | Low |
| 1 | 6.5 | 2.7 | 9.2 | 4.4 | 25.4 | 14.4 |
| 2 | 6.7 | 2.8 | 7.9 | 3.3 | 22.4 | 11.4 |
| 3 | 5.4 | 1.9 | 9.5 | 4.6 | 20.1 | 9.4 |
| 4 | 4.9 | 1.6 | 7.7 | 3.2 | 19 | 8 |
| 5 | 4.8 | 1.5 | 7.3 | 2.8 | 22.1 | 11.1 |
| 6 | 4.4 | 1.3 | 7.5 | 3 | 25.2 | 14.2 |
| 7 | 5.5 | 2 | 6.5 | 2.2 | 18.4 | 7.4 |
| 8 | 6.3 | 2.5 | 9.5 | 4.6 | 20.8 | 9.8 |